# A FUNCTIONAL CENTRAL LIMIT THEOREM FOR THE M/GI/∞ QUEUE

By Laurent Decreusefond and Pascal Moyal

*GET, Télécom Paris, and Université Paris, Dauphine and Université de Technologie de Compiègne*

In this paper, we present a functional fluid limit theorem and a functional central limit theorem for a queue with an infinity of servers M/GI/∞. The system is represented by a point-measure valued process keeping track of the remaining processing times of the customers in service. The convergence in law of a sequence of such processes after rescaling is proved by compactness-uniqueness methods, and the deterministic fluid limit is the solution of an integrated equation in the space $\mathcal{S}'$ of tempered distributions. We then establish the corresponding central limit theorem, that is, the approximation of the normalized error process by a $\mathcal{S}'$-valued diffusion. We apply these results to provide fluid limits and diffusion approximations for some performance processes.

**1. Introduction.** The queues with an infinite reservoir of servers are classical models in queueing theory. In such cases, all the customers are immediately taken care of upon arrival, and spend in the system a sojourn time equal to their service time. Beyond its interest to represent telecommunication networks or computers architectures in which the number of resources is extremely large, this model (commonly referred to as *pure delay* queue) has been often used for comparison to other ones whose dynamics are formally much more complicated, but close in some sense. Then the performances of the pure delay system may give good estimators, or bounds, of that of the other system.

The studies proposed in the literature mainly focused on classical descriptors, such as the length of the queue: among others, its stationary regime under Markovian assumptions ([11]), the transient behavior and law of hitting times of given levels ([14]), the fluid limit and diffusion approximations of normalized sequences ([3, 15]) are now classical results. Let us also mention the recent study on the existence/uniqueness of a stationary regime for the process counting the largest remaining processing time of a customer in the system ([21]).

But when one aim is to have more accurate information on the state of the system, such as the total amount of work at current time (*workload*) or the number of customers having remaining processing time in a given range, no such simple state descriptor can be used. In order to address such questions, one has to know at current time the exhaustive collection of residual processing times of all the customers present in the system. Consequently, we represent the queue by a point measure-valued process $(\mu_t)_{t\geq 0}$, putting Dirac measures at all residual processing times. The price to pay to have such a global information is therefore to work on a very big state space, in fact, of infinite dimension. In the past fifteen years an increasing interest has been dedicated to the study of such measure-valued Markov processes (for definition and main properties, see the reference survey of Dawson [5] on this subject). Such a framework is particularly adequate to describe particles or branching systems (see [5, 22, 24]), or queueing systems whose dynamics are too complex to be carried on with simple finite-dimensional processes: the processor sharing queue (see [12, 26]), queues with deadlines (see [6] for a queue under the earliest deadline first service discipline without reneging, [7, 8] for the same system with reneging and [13] for a processor sharing queue with reneging), or the Shortest Remaining Processing time queue ([1, 20]).

In this paper we aim to identify the "mean behavior" of the measure valued process $(\mu_t)_{t\geq 0}$ describing the pure delay system, introduced above. In that purpose, we use the recent tools of normalization of processes to identify the fluid limit of the process, or formal law of large numbers. This fluid limit is the continuous and deterministic limit in law of a normalized sequence of these processes. We characterize as well the accuracy of this approximation by providing the corresponding functional central limit theorem, that is, the convergence in law of the normalized process of difference between the normalized process and its fluid limit to a diffusion.

Formally, it is rather straightforward in our case to identify the infinitesimal generator of the Feller process $(\mu_t)_{t\geq 0}$. The natural but unusual term is that due to the continuous decreasing of the residual processing times at unit rate as time goes on. This term involves a "spatial derivative" of the measure $\mu_t$, a notion which can only be rigorously defined within the framework of distributions. Because of this term, the fluid limit equation [see (4)] is the integrated version of a partial differential equation rather than an ordinary differential equation, as it is the rule in the previously studied queueing systems. Thus, the classical Gronwall's lemma is of no use here. Fortunately, we can circumvent this difficulty by solving the involved integrated equation, known as transport equation, which is simple enough to have a closed form solution—see Theorem 1. Then, we can proceed using more classical techniques to show the convergence in law of the normalized sequence to the fluid limit (the authors have been informed during the review process of this paper that a similar result has been announced independently in [13]). As a second step, we show the weak convergence of the normalized sequence of deviation to the limit to a diffusion process.

This paper is organized as follows. After some preliminaries in Section 2, we define properly the *profile process* $(\mu_t)_{t\geq 0}$ in Section 3, show, in particular, that $(\mu_t)_{t\geq 0}$ is Feller–Dynkin, and give the corresponding martingale property. In Section 4, we give the fluid limit of $(\mu_t)_{t\geq 0}$. We deduce from this result fluid approximations of some performance processes in Section 5. We prove the functional central limit theorem for $(\mu_t)_{t\geq 0}$ in Section 6, and give the diffusion approximations of the performance processes in Section 6.3.

**2. Preliminaries.** We denote by $\mathcal{D}_b$ (resp. $\mathcal{C}_b$, $\mathcal{C}_K$) the set of real-valued functions defined on $\mathbb{R}$ which are bounded, right-continuous with left-limit (rcll for short) (resp. bounded continuous, continuous with compact support). The space $\mathcal{D}_b$ is equipped with the Skorokhod topology and $\mathcal{C}_b$ with the topology of the uniform convergence. The space of bounded (resp. with compact support) differentiable functions from $\mathbb{R}$ to itself is denoted by $\mathcal{C}_b^1$ (resp. $\mathcal{C}_K^1$) and for $\phi \in \mathcal{C}_b^1$,

$$\|\phi\|_\infty := \sup_{x\in\mathbb{R}}(|\phi(x)| + |\phi'(x)|).$$

We denote by **1** the real function constantly equal to 1, $I$, the identity function $I(x) = x$, $x \in \mathbb{R}$, and for all Borel set $B \in \mathfrak{B}(\mathbb{R})$, $\mathbf{1}_B$ the indicator function of $B$. For all $f \in \mathcal{D}_b$ and all $x \in \mathbb{R}$, we denote by $\tau_x f$, the function $\tau_x f(\cdot) := f(\cdot - x)$.

The Schwartz space, denoted by $\mathcal{S}$, is the space of infinitely differentiable functions, equipped with the topology defined by the semi-norms

$$|\phi|_{\beta,\gamma} := \sup_{x\in\mathbb{R}}\left|x^\beta \frac{d^\gamma}{dx^\gamma}\phi(x)\right|, \qquad \beta \in \mathbb{N}, \gamma \in \mathbb{N}.$$

Its topological dual, the space of tempered distributions, is denoted by $\mathcal{S}'$, and the duality product is classically denoted $\langle\mu,\phi\rangle$. The *distributional derivative* of some $\mu \in \mathcal{S}'$ is $\mu' \in \mathcal{S}'$ such that $\langle\mu',\phi\rangle = -\langle\mu,\phi'\rangle$ for all $\phi \in \mathcal{S}$. For all $\mu \in \mathcal{S}'$ and $x \in \mathbb{R}$, we denote by $\tau_x\mu$ the tempered distribution satisfying $\langle\tau_x\mu,\phi\rangle = \langle\mu,\tau_x\phi\rangle$ for all $\phi \in \mathcal{S}$.

The set of finite nonnegative measures on $\mathbb{R}$ is denoted by $\mathcal{M}_f^+$ (which is part of $\mathcal{S}'$) and $\mathcal{M}_p$ is the set of finite counting measures on $\mathbb{R}$. The space $\mathcal{M}_f^+$ is equipped with the weak topology $\sigma(\mathcal{M}_f^+, \mathcal{C}_b)$, for which $\mathcal{M}_f^+$ is Polish (we write $\langle \mu, f \rangle = \int f \, d\mu$ for $\mu \in \mathcal{M}_f^+$ and $f \in \mathcal{D}_b$). We say that a sequence $\{\mu^n\}_{n \in \mathbb{N}^*}$ of $\mathcal{M}_f^+$ converges *weakly* to $\mu$, and denote $\mu^n \overset{w}{\Rightarrow} \mu$, if for all $f \in \mathcal{C}_b$, $\langle \mu^n, f \rangle$ tends to $\langle \mu, f \rangle$.

We also denote for all $x \in \mathbb{R}$ and all $\nu \in \mathcal{M}_f^+$, $\tau_x \nu$ the measure satisfying for all Borel set $B$, $\tau_x \nu(B) := \nu(B - x)$. Remark that these last two definitions coincide with that in $\mathcal{S}'$: for all $\phi \in \mathcal{S}'$ and $\nu \in \mathcal{M}_f^+$, $\langle \nu, \phi \rangle_{\mathcal{S}'} = \langle \nu, \phi \rangle_{\mathcal{M}_f^+}$, and for all $x \in \mathbb{R}$, $\langle \tau_x \nu, \phi \rangle = \langle \nu, \tau_x \phi \rangle$.

For a random variable (r.v. for short) $X$ defined on a fixed probability space $(\Omega, \mathcal{F})$, we say that $\{X^n\}_{n \in \mathbb{N}^*}$ converges *in distribution* to $X$, and denote $X^n \Rightarrow X$, if the sequence of the distributions of the $X^n$'s tends weakly to that of $X$.

Let $\mathcal{C}(\mathcal{M}_f^+, \mathbb{R})$ be the set of continuous functions from $\mathcal{M}_f^+$ to $\mathbb{R}$. Letting $0 < T \le \infty$, for $E$ a Polish space, we denote $\mathcal{C}([0, T], E)$, respectively $\mathcal{D}([0, T], E)$, the Polish space (for its usual strong topology) of continuous, respectively rcll, functions from $[0, T]$ to $E$. The mutual variation of two local martingales $(M_t)_{t \ge 0}$ and $(N_t)_{t \ge 0}$ in $\mathcal{D}([0, T], \mathbb{R})$ is denoted by $(<M, N>_t)_{t \ge 0}$ with $(<M>_t)_{t \ge 0} = (<M, M>_t)_{t \ge 0}$ for the quadratic variation (or the increasing process) of $(M_t)_{t \ge 0}$.

Since $\mathcal{S}'$ is a nuclear Fréchet space (see [27]), note that $(\mu_t)_{t \ge 0} \in \mathcal{C}([0, T], \mathcal{S}')$ (resp. $\mathcal{D}([0, T], \mathcal{S}')$) if and only if $(\langle \mu_t, \phi \rangle)_{t \ge 0} \in \mathcal{C}([0, T], \mathbb{R})$ (resp. $\mathcal{D}([0, T], \mathbb{R})$) for all $\phi \in \mathcal{S}$.

Let $T > 0$, and let $(X_t)_{t \ge 0}$[1] be a process of $\mathcal{D}([0, T], \mathcal{S}')$. We denote for all $t \in [0, T]$, $\int_0^t X_s \, ds$, the element of $\mathcal{S}'$ such that, for all $\phi \in \mathcal{S}$,

$$\left\langle \int_0^t X_s \, ds, \phi \right\rangle = \int_0^t \langle X_s, \phi \rangle \, ds.$$

Let us, moreover, define the following $\mathcal{S}'$-valued processes:

(1) $\qquad (G(X)_t)_{t \ge 0} = (\tau_t X_t)_{t \ge 0},$

(2) $\qquad (H(X)_t)_{t \ge 0} = \left( \int_0^t X_s \, ds \right)_{t \ge 0},$

(3) $\qquad (N(X)_t)_{t \ge 0} = \left( \int_0^t \tau_{t-s}(X_s)' \, ds \right)_{t \ge 0},$

which means, for all $\phi \in \mathcal{S}$ and all $t \in [0, T]$,

$$\langle N(X)_t, \phi \rangle = -\int_0^t \langle X_s, \tau_{t-s} \phi' \rangle \, ds.$$

Let $(\mathcal{F}_t)_{t \ge 0}$ be a filtration on $(\mathbb{R}, \mathfrak{B}(\mathbb{R}))$. Let us recall (cf. [17, 27]) that a process $M \in \mathcal{D}([0, T], \mathcal{S}')$ is a $\mathcal{S}'$ valued $\mathcal{F}_t$-*semi-martingale* (*resp. local martingale, martingale*) if, for all $\phi \in \mathcal{S}$, $(\langle M_t, \phi \rangle)_{t \ge 0}$ is a real $\mathcal{F}_t$-semi-martingale (resp. local martingale, martingale). According to [18], page 13, the *tensor-quadratic* process $<<M>>$ of the $\mathcal{S}'$ valued martingale $M$ is given for all $t \geq 0$ and all $\phi, \psi \in \mathcal{S}$ by

$$<<M>>_t(\phi, \psi) := <\langle M_\cdot, \phi\rangle, \langle M_\cdot, \psi\rangle>_t.$$

---

[1]when no ambiguity is possible, we often write for instance $X$ for the process $(X_t)_{t \ge 0}$.

For all $\mathcal{S}'$-valued semi-martingale $M$, and all real semi-martingale $X$, the $\mathcal{S}'$-valued *stochastic integral of $M$ with respect to $X$* is denoted by $\int_0^t M_s\, dX_s$, and is such that, for all $t \geq 0$ and all $\phi \in \mathcal{S}$,

$$\left\langle \int_0^t M_s\, dX_s, \phi \right\rangle = \int_0^t \langle M_s, \phi\rangle\, dX_s.$$

Let $T \geq 0$. We say that $X \in \mathcal{D}([0,T], \mathcal{S}')$ satisfies the *integrated transport equation* associated to $(K, g)$ (see [8]) if, for some $K \in \mathcal{S}'$ and $(g_t)_{t \geq 0} \in \mathcal{D}([0,T], \mathcal{S}')$,

(4) $$X_t = K + \int_0^t (X_s)'\, ds + g_t \qquad \text{for all } t \in [0,T].$$

Let us then recall the following result.

THEOREM 1 ([8], Theorem 1). *The only solution in $\mathcal{D}([0,T], \mathcal{S}')$ of the integrated transport equation associated to $(K, g)$ is given for all $t \in [0, T]$ by*

(5) $$X_t = \tau_t K + g_t + N(g)_t,$$

*where the mapping $N$ is defined in (3).*

Note that Propositions 3 and 4 show, in particular, that $X$ defined by (5) is an element of $\mathcal{D}([0,T], \mathcal{S}')$.

**3. The profile process of the M/GI/∞ queue.** Hereafter we consider a pure delay system M/GI/∞: on a probability space $(\Omega, \mathcal{F}, \mathbf{P})$, consider a Poisson process

$$N_t := \sum_{i \in \mathbb{N}^*} \mathbf{1}_{\{T_i \leq t\}}$$

of positive intensity $\lambda$, representing the arrivals of the customers in a queueing system with an infinite reservoir of servers. Hence, any of them is immediately attended upon arrival. The time spent in the system by the $i$th arriving customer (denoted $C_i$) equals the service duration $\sigma_i$ he requests. We assume that the sequence of marks $\{\sigma_i\}_{i \in \mathbb{N}^*}$ is i.i.d. with the random variable $\sigma$ having the distribution $\alpha \in \mathcal{M}_f^+$. We denote at all time $t \geq 0$, $X_t$ the number of customers currently in the system, and $S_t$ the number of already served customers at $t$, related to $N_t$ and $X_t$ by the relation $N_t = X_t + S_t$. The *workload process* $(W_t)_{t \geq 0}$ equals at time $t$ the total amount of service requested by the customers in the system at $t$, in time units. The *profile process* of the queue is the point-measure valued process $(\mu_t)_{t\geq 0}$ whose units of mass represent the remaining processing times (i.e., time to the service completion) of all the customers who already entered the system at current time. In other words, for all $t \geq 0$,

$$\mu_t = \sum_{i=1}^{N_t} \delta_{(T_i+\sigma_i-t)}.$$

In this expression a nonpositive remaining processing time $(T_i + \sigma_i - t \leq 0)$ stands for a customer who already left the system, precisely $t - (T_i + \sigma_i)$ time units before $t$. As easily seen, for all $t \geq 0$ and any $x < y$, the number of customers having at $t$ a remaining processing time in $(x,y)$ (the boundaries may be included or infinite) is given by $\langle \mu_t, \mathbf{1}_{(x,y)} \rangle$. Thus, the congestion, service and workload processes $X$, $S$ and $W$ are easily recovered from $\mu$ by the relations $X_t = \langle \mu_t, \mathbf{1}_{\mathbb{R}_+^*} \rangle$, $S_t = \langle \mu_t, \mathbf{1}_{\mathbb{R}_-} \rangle$ and $W_t = \langle \mu_t, I\mathbf{1}_{\mathbb{R}_+^*} \rangle$, $t \geq 0$.

The dynamics of $(\mu_t)_{t\geq 0}$ can be described as follows. For any $t \geq 0$, $\mu_{t+dt}$ is the translated of $\mu_t$ toward left by $dt$, that is, $\tau_{dt}\mu_t$ (the remaining service times of all the customers in the system at $t$, if any, decrease at unit rate), and a Dirac measure at $\sigma$ is added if a new customer enters the system between $t$ and $t + dt$ having service time $\sigma$. In other words, we have for any $\phi \in \mathcal{C}_K$

(6) $$\langle \mu_{t+dt}, \phi \rangle = \langle \tau_{dt}\mu_t, \phi \rangle + \phi(\sigma) N(dt, d\sigma),$$

where $N(\cdot, \cdot)$ denotes the Poisson measure associated to the marked arrival point process. Hence, the above equation allows one to construct $(\langle \mu_t, \phi \rangle)_{t\geq 0}$ for any $\phi \in \mathcal{C}_K$ by induction on the arrival times and for any initial state $\mu \in \mathcal{M}_f^+$ (as in [26], page 192). Since $\mathcal{C}_K$ is a separating class of $\mathcal{M}_f^+$, the initial value $\mu_0$ and the evolution equation (6) fully define $(\mu_t)_{t\geq 0}$. Let $h > 0$ and $F$ be a bounded continuous function from $\mathcal{M}_f^+$ into $\mathbb{R}$. Denote for all $i \in \mathbb{N}$

$$\mathcal{A}_i(h) = \{\text{There are exactly } i \text{ arrivals in } [0,h]\}.$$

We have for all $\mu \in \mathcal{M}_f^+$

$$T_h F(\mu) := \mathbf{E}[F(\mu_h) | \mu_0 = \mu]$$

(7) $$= \sum_{i\geq 0} \mathbf{E}\left[ F\left(\tau_h \mu + \sum_{j>0} \delta_{\sigma_j - (h-T_j)} \mathbf{1}_{\{T_j \leq h\}}\right) \mathbf{1}_{\mathcal{A}_i(h)} \,\middle|\, \mu_0 = \mu \right]$$

$$= (1 - \lambda h) F(\tau_h \mu) + \lambda h \int F(\tau_h(\mu + \delta_x)) \, d\alpha(x) + \epsilon(\mu, h),$$

where

$$
\begin{aligned}
\epsilon(\mu, h) = &\{\mathbf{P}[\mathcal{A}_0(h)] - (1 - \lambda h)\} F(\tau_h \mu) \\
&+ \{\mathbf{P}[\mathcal{A}_1(h)] - \lambda h\} \int F(\tau_h(\mu + \delta_x))\, d\alpha(x) \\
&+ \mathbf{E}[\{F(\tau_h \mu + \delta_{\sigma_1 - (h - T_1)}) - F(\tau_h(\mu + \delta_{\sigma_1}))\} \mathbf{1}_{\mathcal{A}_1(h)} \mid \mu_0 = \mu] \\
&+ \sum_{i \geq 2} \mathbf{E}\left[ F\left(\tau_h \mu + \sum_{j=1}^{i} \delta_{\sigma_j - (h - T_j)}\right) \mathbf{1}_{\mathcal{A}_i(h)} \mid \mu_0 = \mu \right].
\end{aligned}
$$

(8)

Hence, $\epsilon(\mu, h)$ is a $o(h)$, as easily seen using dominated convergence for the third term on the right-hand side of (8), and from the fact that $F$ is bounded for the other three terms. Therefore, according to [5], page 18, $(\mu_t)_{t \geq 0}$ is a weak homogeneous $\mathcal{M}_f^+$-valued Markov process having transition $T_h F$ for all $h > 0$ and all bounded continuous $F$. Note, moreover, that $(\mu_t)_{t \geq 0} \in \mathcal{D}([0, \infty), \mathcal{M}_f^+)$ since $(\langle \mu_t, \phi \rangle)_{t \geq 0} \in \mathcal{D}([0, \infty), \mathbb{R})$ for all $\phi \in \mathcal{C}_b$.

PROPOSITION 1. *The process $(\mu_t)_{t \geq 0}$ is a Feller–Dynkin process of $\mathcal{D}([0, \infty), \mathcal{M}_f^+)$.*

PROOF. According to Lemma 3.5.1 and Corollary 3.5.2 of [5], the Markov process $(\mu_t)_{t \geq 0}$ enjoys the Feller–Dynkin property if:

(i) For all $f \in \mathcal{C}_K^1$, $h > 0$, the mapping $\mu \mapsto \mathbf{E}[F_f(\mu_h) \mid \mu_0 = \mu]$ is continuous.
(ii) For all $h > 0$, $\mathbf{E}[F_{\mathbf{1}}(\mu_h) \mid \mu_0 = \mu] \longrightarrow 0$, as $\mu(\mathbb{R}) \to +\infty$.
(iii) For all $\mu \in \mathcal{M}_f^+$ and $f \in \mathcal{C}_K^1$, $\mathbf{E}[F_f(\mu_h) \mid \mu_0 = \mu] \longrightarrow F_f(\mu)$ as $h \to 0$.

Let us denote for all $f \in \mathcal{C}_K^1$, $F_f$ the mapping from $\mathcal{M}_f^+$ into $\mathbb{R}$ defined by $F_f(\mu) = e^{-\langle \mu, f \rangle}$. In view of (7) and (8), we have for all $\mu \in \mathcal{M}_f^+$ and $f \in \mathcal{C}_K^1$

(9)
$$
\begin{aligned}
&\mathbf{E}[F_f(\mu_h) \mid \mu_0 = \mu] \\
&= e^{-\langle \tau_h \mu, f \rangle}\left(1 - \lambda h + \lambda h \int e^{-\langle \tau_h \delta_x, f \rangle}\, d\alpha(x) + \epsilon_f(h)\right),
\end{aligned}
$$

where

$$
\begin{aligned}
\epsilon_f(h) = &\{\mathbf{P}[\mathcal{A}_0(h)] - (1 - \lambda h)\} + \{\mathbf{P}[\mathcal{A}_1(h)] - \lambda h\} \int e^{-f(x-h)}\, d\alpha(x) \\
&+ \mathbf{E}[\{e^{-f(\sigma_1 - (h - T_1))} - e^{-f(\sigma_1 - h)}\} \mathbf{1}_{\mathcal{A}_1(h)}] \\
&+ \sum_{i \geq 2} \mathbf{E}\left[\left\{\prod_{j=1}^{i} e^{-f(\sigma_j - (h - T_j))}\right\} \mathbf{1}_{\mathcal{A}_i(h)}\right],
\end{aligned}
$$

which is a $o(h)$ from the same arguments as for (8). Note, moreover, that $\epsilon_f(h)$ does not depend on $\mu$. Hence, the continuity in (i) is granted by the fact that, for all $h > 0$, the mapping $\mu \mapsto \tau_h \mu$ is continuous from $\mathcal{M}_f^+$ into itself. For (ii), remark that the total mass on $\mathbb{R}$ of $\tau_h \mu$ equals that of $\mu$. Finally, the convergence in (iii) holds since $\langle \mu, \tau_h f \rangle \longrightarrow_{h \to 0} \langle \mu, f \rangle$ by dominated convergence. The proposition is proved. $\square$

As easily seen from (9), and since the set of linear combinations of $F_f$, $f \in \mathcal{C}_K^1$, is dense in $\mathcal{C}(\mathcal{M}_f^+, \mathbb{R})$ (see [26], Proposition 7.10), the infinitesimal generator $\mathcal{A}$ of $(\mu_t)_{t \geq 0}$ is given for all $\mu \in \mathcal{M}_f^+$ by

$$\mathcal{A}F(\mu) := \lim_{h \to 0} \frac{\mathbf{E}[F(\mu_h) \mid \mu_0 = \mu] - F(\mu)}{h}$$
$$= \lim_{h \to 0} \frac{F(\tau_h \mu) - F(\mu)}{h} - \lambda F(\mu) + \lambda \int F(\mu + \delta_x) \, d\alpha(x),$$

for all $F \in \mathcal{C}(\mathcal{M}_f^+, \mathbb{R})$ such that the latter limit exists (we say that $F$ belongs to the *domain* of $\mathcal{A}$).

PROPOSITION 2. *For all $\phi \in \mathcal{C}_b^1$, the process defined for all $t \geq 0$ by*

$$M_t(\phi) = \langle \mu_t, \phi \rangle - \langle \mu_0, \phi \rangle + \int_0^t \langle \mu_s, \phi' \rangle \, ds - \lambda t \langle \alpha, \phi \rangle$$

*is an rcll square integrable $\mathcal{F}_t$-martingale. For all $\phi, \psi \in \mathcal{S}$, the mutual variation of $(M_t(\phi))_{t \geq 0}$ with $(M_t(\psi))_{t \geq 0}$ is given for all $t \geq 0$ by*

$$< M_\cdot(\phi), M_\cdot(\psi) >_t = \lambda t \langle \alpha, \phi \psi \rangle. \tag{10}$$

PROOF. For all $\phi \in \mathcal{C}_b^1$, the mapping $\Pi_\phi$ from $\mathcal{M}_f^+$ into $\mathbb{R}$ defined by $\Pi_\phi(\mu) = \langle \mu, \phi \rangle$ clearly belongs to the domain of $\mathcal{A}$. As a consequence of Dynkin's lemma ([9, 10]), the process defined for all $t \geq 0$ by

$$M_t(\phi) = \Pi_\phi(\mu_t) - \Pi_\phi(\mu_0) - \int_0^t \mathcal{A}\Pi_\phi(\mu_s) \, ds$$
$$= \langle \mu_t, \phi \rangle - \langle \mu_0, \phi \rangle + \int_0^t \langle \mu_s, \phi' \rangle \, ds - \lambda t \langle \alpha, \phi \rangle \tag{11}$$

is an rcll $\mathcal{F}_t$-local martingale. Let now $\psi \in \mathcal{C}_b^1$. The mapping $\Pi_{\phi,\psi}$ from $\mathcal{M}_f^+$ into $\mathbb{R}$ defined by $\Pi_{\phi,\psi}(\mu) = \Pi_\phi(\mu) \Pi_\psi(\mu)$ also belongs to the domain of $\mathcal{A}$, implying that

$$\tilde{M}_t(\phi, \psi) = \Pi_\phi(\mu_t) \Pi_\psi(\mu_t) - \Pi_\phi(\mu_0) \Pi_\psi(\mu_0) - \int_0^t \mathcal{A}\Pi_{\phi,\psi}(\mu_s) \, ds \tag{12}$$

is as well an rcll $\mathcal{F}_t$-local martingale. But as easily checked, for all $\mu \in \mathcal{M}_f^+$,

(13) $\quad \mathcal{A}\Pi_{\phi,\psi}(\mu) = \Pi_\phi(\mu)\mathcal{A}\Pi_\psi(\mu) + \Pi_\psi(\mu)\mathcal{A}\Pi_\phi(\mu) + \lambda\langle\alpha,\phi\psi\rangle.$

Itô's formula yields together with (11) that, for all $t \geq 0$,

$\Pi_\phi(\mu_t)\Pi_\psi(\mu_t)$

$$= \Pi_\phi(\mu_0)\Pi_\psi(\mu_0) + \int_0^t \Pi_\phi(\mu_s)\,dM_s(\psi) + \int_0^t \Pi_\phi(\mu_s)\mathcal{A}\Pi_\psi(\mu_s)\,ds$$

$$+ \int_0^t \Pi_\psi(\mu_s)\,dM_s(\phi) + \int_0^t \Pi_\psi(\mu_s)\mathcal{A}\Pi_\phi(\mu_s)\,ds$$

$$+ <\Pi_\phi(\mu_\cdot), \Pi_\psi(\mu_\cdot)>_t,$$

which, together with (12) and (13), implies

$$\int_0^t \Pi_\phi(\mu_s)\,dM_s(\psi) + \int_0^t \Pi_\psi(\mu_s)\,dM_s(\phi) + <\Pi_\phi(\mu_\cdot), \Pi_\psi(\mu_\cdot)>_t$$

$$= \tilde{M}_t(\phi,\psi) + \lambda t\langle\alpha,\phi\psi\rangle.$$

By identifying the finite variation processes, we obtain that **P**-a.s. for all $t \geq 0$,

(14) $\quad <\Pi_\phi(\mu_\cdot), \Pi_\psi(\mu_\cdot)>_t = \lambda t\langle\alpha,\phi\psi\rangle,$

but in view of (11), this last quantity equals $<M_\cdot(\phi), M_\cdot(\psi)>_t$. In particular, for all $t \geq 0$,

$$\mathbf{E}[<M_\cdot(\phi)>_t] = \lambda t\langle\alpha,\phi^2\rangle < \infty,$$

hence, $(M_t(\phi))_{t\geq 0}$ is a square integrable martingale. $\square$

**4. Fluid limit.** Consider a sequence of M/GI/$\infty$ systems such that the $n$th system has an initial profile $\mu_0^n$, is fed by a Poisson process $(N_t^n)_{t\geq 0}$ of arrival times $\{T_i^n\}_{i\in\mathbb{N}^*}$ and of intensity $\lambda^n$, in which the customers request i.i.d. service durations $\{\sigma_i^n\}_{i\in\mathbb{N}^*}$ which have the nonatomic distribution $\alpha^n$ of a r.v. $\sigma^n$. We assume, furthermore, that $\sigma^n$ is integrable, that is, $\langle\alpha^n, I\rangle < \infty$. The process $(\mu_t^n)_{t\geq 0}$, defined for all $t$ by

$$\mu_t^n = \sum_{i=1}^{N_t^n} \delta_{T_i^n + \sigma_i^n - t},$$

is the profile process of this $n$th system. Denote by $(\mathcal{F}_t^n)_{t\geq 0}$ the associated filtration. Let us also define as previously the performance processes of the $n$th system: $X^n$, $S^n$ and $W^n$. We normalize the process $\mu^n$ in time, space and weight by defining for all $t \geq 0$

$$\bar{\mu}_t^n = \frac{1}{n}\sum_{i=1}^{N_{nt}^n} \delta_{(T_i^n + \sigma_i^n - nt)/n}.$$

Thus, for all Borel set $B$ and all $t$,

$$\bar{\mu}_t^n(B) = \frac{\mu_{nt}^n(nB)}{n},$$

where $nB := \{nx, x \in B\}$. In words, $\bar{\mu}_t^n$ is $1/n$ times the point measure having atoms at levels $n$ times smaller than that of $\mu_{nt}^n$. We also define $(\mathcal{G}_t^n)_{t \geq 0} := (\mathcal{F}_{nt}^n)_{t \geq 0}$, the associated filtration, and normalize the arrival process as well as the performance processes of the $n$th system the corresponding way, that is, for all $t \geq 0$,

$$\bar{N}_t^n := \frac{N_{nt}^n}{n}, \qquad \bar{X}_t^n := \frac{X_{nt}^n}{n} = \langle \bar{\mu}_t^n, \mathbf{1}_{\mathbb{R}_+^*} \rangle,$$

$$\bar{S}_t^n := \frac{S_{nt}^n}{n} = \langle \bar{\mu}_t^n, \mathbf{1}_{\mathbb{R}_-} \rangle, \qquad \bar{W}_t^n := \frac{W_{nt}^n}{n^2} = \langle \bar{\mu}_t^n, I\mathbf{1}_{\mathbb{R}_+^*} \rangle.$$

We denote for all $i \in \mathbb{N}^*$, $\bar{\sigma}_i^n := \sigma_i^n/n$. Thus, the $\{\bar{\sigma}_i^n\}_{i \in \mathbb{N}^*}$ are i.i.d. with the nonatomic distribution $\bar{\alpha}^n$ of the r.v. $\sigma^n/n$. The distribution $\bar{\alpha}^n$ is such that $\langle \bar{\alpha}^n, I \rangle < \infty$, and satisfies for all Borel set $B$

$$\bar{\alpha}^n(B) = \alpha^n(nB).$$

We assume hereafter that the following hypothesis holds.

HYPOTHESIS 1. There exists $\lambda > 0$ such that

(15) $$\lambda^n \xrightarrow[n \to \infty]{} \lambda,$$

For all $\varepsilon > 0$, there exists $M_\varepsilon > 0$ such that, for all $n \in \mathbb{N}^*$,

(16) $$\mathbf{P}[\langle \mu_0^n, \mathbf{1} \rangle > nM_\varepsilon] \leq \varepsilon.$$

For some measure $\bar{\mu}_0^*$ of $\mathcal{M}_f^+$ such that

$$\langle \bar{\mu}_0^*, I \rangle < \infty,$$

for all $f \in \mathcal{D}_b$,

(17) $$\langle \bar{\mu}_0^n, f \rangle \xrightarrow[n \to \infty]{} \langle \bar{\mu}_0^*, f \rangle \qquad \text{in probability.}$$

There exists a nonatomic probability distribution $\bar{\alpha}^*$ such that

$$\bar{\alpha}^n \xRightarrow{w} \bar{\alpha}^*,$$

$$\langle \bar{\alpha}^n, I \rangle \xrightarrow[n \to \infty]{} \langle \bar{\alpha}^*, I \rangle < \infty.$$

Let $\phi \in \mathcal{C}_b^1$ and $\psi^n(\cdot) = \phi(\cdot/n)/n$. According to Proposition 2, the process defined for all $t$ by

(18) $$\bar{M}_t^n(\phi) := M_{nt}^n(\psi^n) = \langle \bar{\mu}_t^n, \phi \rangle - \langle \bar{\mu}_0^n, \phi \rangle + \int_0^t \langle \bar{\mu}_s^n, \phi' \rangle \, ds - \lambda^n t \langle \bar{\alpha}^n, \phi \rangle$$

is a square integrable $\mathcal{G}_t^n$-martingale of $\mathcal{D}([0,\infty),\mathbb{R})$, such that

(19) $$< \bar{M}_\cdot^n(\phi) >_t = \frac{\lambda^n}{n} t \langle \bar{\alpha}^n, \phi^2 \rangle.$$

In other words, the process $\bar{M}^n$ defined for all $t$ by

$$\bar{M}_t^n = \bar{\mu}_t^n - \bar{\mu}_0^n - \int_0^t (\bar{\mu}_s^n)' \, ds - \lambda^n t \bar{\alpha}^n$$

is a $\mathcal{S}'$-valued $\mathcal{G}_t^n$-martingale of tensor-quadratic process given for any $\phi$ and $\psi$ in $\mathcal{S}$ by

$$<< \bar{M}^n >>_t (\phi,\psi) = \frac{\lambda^n}{n} t \langle \bar{\alpha}^n, \phi\psi \rangle.$$

THEOREM 2. *Assume that Hypothesis 1 holds. Then*
$$\bar{\mu}^n \Longrightarrow \bar{\mu}^* \qquad \text{in } \mathcal{D}([0,\infty), \mathcal{M}_f^+),$$

*where $\bar{\mu}^*$ is the deterministic element of $\mathcal{C}([0,\infty), \mathcal{M}_f^+)$ defined for all $t \geq 0$ and all $\phi \in \mathcal{D}_b$ by*

(20) $$\langle \bar{\mu}_t^*, \phi \rangle = \langle \bar{\mu}_0^*, \tau_t \phi \rangle + \lambda \int_0^t \langle \bar{\alpha}^*, \tau_s \phi \rangle \, ds.$$

PROOF. First we prove that $\{\bar{\mu}^n\}_{n \in \mathbb{N}^*}$ is tight in $\mathcal{D}([0,\infty), \mathcal{M}_f^+)$. To this end, it suffices to show conditions **C.1** and **C.2** of Theorem A.2 in the Appendix. To show **C.1**, fix $\phi \in \mathcal{C}_b^1$ and $T > 0$. Remarking that, for all $s \geq 0$,

$$\langle \bar{\mu}_s^n, \mathbf{1} \rangle \leq \bar{N}_s^n + \langle \bar{\mu}_0^n, \mathbf{1} \rangle,$$

equation (18) yields **P**-a.s. for all $u < v \leq T$,

(21) $$\begin{aligned}|\langle \bar{\mu}_v^n, \phi \rangle - \langle \bar{\mu}_u^n, \phi \rangle| \\ &\leq \int_u^v |\langle \bar{\mu}_s^n, \phi' \rangle| \, ds + \lambda^n |\langle \bar{\alpha}^n, \phi \rangle||v-u| + |\bar{M}_\phi^n(v) - \bar{M}_\phi^n(u)| \\ &\leq |v-u|\{\|\phi'\|_\infty (\bar{N}_T^n + \langle \bar{\mu}_0^n, \mathbf{1} \rangle) + \|\phi\|_\infty \lambda^n\} \\ &\quad + |\bar{M}_\phi^n(v) - \bar{M}_\phi^n(u)|.\end{aligned}$$

Let $\xi > 0$. From Doob's inequality,

$$\mathbf{P}\left[\sup_{t \leq T} |\bar{M}_t^n(\phi)| \geq \xi\right] \leq \frac{4}{\xi^2} \mathbf{E}[< \bar{M}_\cdot^n(\phi) >_T] \leq \frac{4}{\xi^2}\left(\frac{\lambda^n}{n}\right) \|\phi^2\|_\infty T \xrightarrow[n \to \infty]{} 0,$$

in view of (19). Hence, it follows from the standard convergence criterion on $\mathcal{D}([0,T],\mathbb{R})$ that $\{(\bar{M}_t^n(\phi))_{t \geq 0}\}_{n \in \mathbb{N}^*}$ converges in distribution to the null process. This sequence is, in particular, tight, and it is routine in view of (21) and Hypothesis 1 to check the standard tightness criterion in $\mathcal{D}([0,T],\mathbb{R})$: for all $\varepsilon > 0$ and $\eta > 0$, there exists $\delta > 0$ and $N \in \mathbb{N}$ such that, for all $n \geq N$,

$$(22) \qquad \mathbf{P}\left[\sup_{u,v \leq T, |v-u| \leq \delta} |\langle \bar{\mu}_v^n, \phi \rangle - \langle \bar{\mu}_u^n, \phi \rangle| \geq \eta\right] \leq \varepsilon.$$

Letting $\beta_\varepsilon := M_\varepsilon \|\phi\|_\infty$, assumption (16) implies that, for all $n \in \mathbb{N}^*$,

$$(23) \qquad \mathbf{P}[|\langle \bar{\mu}_0^n, \phi \rangle| > \beta_\varepsilon] \leq \mathbf{P}[\|\phi\|_\infty \langle \bar{\mu}_0^n, \mathbf{1}\rangle > M_\varepsilon \|\phi\|_\infty] \leq \varepsilon.$$

Hence, (22) and (23) show that $\{(\langle \bar{\mu}_t^n, \phi\rangle)_{t \geq 0}\}_{n \in \mathbb{N}^*}$ is tight in $\mathcal{D}([0,T],\mathbb{R})$ for all $T > 0$ (see [25], Theorem D.9). This sequence is thus tight in $\mathcal{D}([0,\infty),\mathbb{R})$.

We now prove condition **C.2** (compact containment). Fix $T > 0$. Let us first apply [12], Lemma A.2.: under Hypothesis 1, we have

$$(24) \qquad \left\{\left(\frac{1}{n}\sum_{i=1}^{n\bar{N}_t^n} \phi(\bar{\sigma}_i^n)\right)_{t \geq 0}\right\}_{n \in \mathbb{N}^*} \Longrightarrow (\lambda t \langle \bar{\alpha}^*, \phi\rangle)_{t \geq 0} \qquad \text{in } \mathcal{D}([0,T],\mathbb{R})$$

for any measurable $\phi$, such that $\phi$ is continuous on the supports of $\bar{\alpha}^*$ and $\bar{\alpha}^n$, $n \in \mathbb{N}^*$, and such that $\langle |\phi|, \bar{\alpha}^*\rangle < \infty$ and $\langle |\phi|, \bar{\alpha}^n\rangle < \infty$, $n \in \mathbb{N}^*$. In particular, this yields for any $0 < l \leq T$, and any such $\phi$,

$$(25) \qquad \mathbf{P}\left[\sup_{t \in [0,T-l]} \frac{1}{n} \sum_{N_{nt}^n+1}^{N_{n(t+l)}^n} \phi(\bar{\sigma}_i^n) > 2\lambda l \langle \bar{\alpha}^*, \phi\rangle\right] \xrightarrow[n \to \infty]{} 0.$$

Taking $l = T$ and $\phi = \mathbf{1}$ (resp. $\phi = I$) in the above expression yields

$$\mathbf{P}[\bar{N}_T^n > 2\lambda T] \xrightarrow[n \to \infty]{} 0,$$

$$\mathbf{P}\left[\frac{1}{n}\sum_{i=1}^{N_{nT}^n} \bar{\sigma}_i^n > 2\lambda T \langle \bar{\alpha}^*, I\rangle\right] \xrightarrow[n \to \infty]{} 0.$$

Letting $M_T = \max\{2\lambda T + \langle \bar{\mu}_0^*, \mathbf{1}\rangle, 2\lambda T \langle \bar{\alpha}^*, I\rangle + \langle \bar{\mu}_0^*, I\rangle\} + 1$, we have

$$\mathbf{P}\left[\sup_{t \in [0,T]} \max\{\langle \bar{\mu}_t^n, \mathbf{1}_{\mathbb{R}_+}\rangle, \langle \bar{\mu}_t^n, I\mathbf{1}_{\mathbb{R}_+}\rangle\} > M_T\right]$$
$$(26) \qquad \leq \mathbf{P}[\langle \bar{\mu}_0^n, \mathbf{1}\rangle > \langle \bar{\mu}_0^*, \mathbf{1}\rangle + 1] + \mathbf{P}[\bar{N}_T^n > 2\lambda T]$$
$$\qquad + \mathbf{P}[\langle \bar{\mu}_0^n, I\rangle > \langle \bar{\mu}_0^*, I\rangle + 1] + \mathbf{P}\left[\frac{1}{n}\sum_{i=1}^{N_{nT}^n} \bar{\sigma}_i^n > 2\lambda T\langle \bar{\alpha}^*, I\rangle\right] \xrightarrow[n \to +\infty]{} 0.$$

For all $0 < \eta < 1$, denote

$$\mathcal{K}_{T,\eta} := \{\zeta \in \mathcal{M}_f^+; \max\{\langle \zeta, \mathbf{1}_{\mathbb{R}_+}\rangle, \langle \zeta, I\mathbf{1}_{\mathbb{R}_+}\rangle\} \leq M_T, \langle \zeta, \mathbf{1}_{(-\infty, -T]}\rangle = 0\}.$$

Since $\langle \zeta, I\mathbf{1}_{\mathbb{R}+}\rangle \leq M_T$ implies that for all $y > 0$, $\zeta([y,\infty)) \leq M_T/y$, we have

$$\lim_{y\to\infty} \sup_{\zeta\in\mathcal{K}_{T,\eta}} \zeta([y,\infty)) = 0, \qquad \lim_{y\to-\infty} \sup_{\zeta\in\mathcal{K}_{T,\eta}} \zeta((-\infty,y]) = 0,$$

which implies in turns that $\mathcal{K}_{T\eta} \subset \mathcal{M}_f^+$ is relatively compact (see [16]). Now, since up to time $T$ no already served customer can have a remaining processing time less than $-T$, we have

$$\sup_{t\leq T}\langle \bar{\nu}_t^n, \mathbf{1}_{(-\infty,-T]}\rangle = 0.$$

This, together with (26), implies that

$$\liminf_{n\to\infty} \mathbf{P}[\bar{\mu}_t^n \in \mathcal{K}_{T,\eta}, \text{ for all } t \in [0,T]] \geq 1 - \eta.$$

$\mathbf{K}_{T,\eta}$ being the closure of $\mathcal{K}_{T,\eta}$, we found a compact subset $\mathbf{K}_{T,\eta} \subset \mathcal{M}_f^+$ such that

$$\liminf_{n\to\infty} \mathbf{P}[\bar{\mu}_t^n \in \mathbf{K}_{T,\eta}, \text{ for all } t \in [0,T]] \geq 1 - \eta,$$

which completes the proof of tightness.

Let now $(\chi_t)_{t\geq 0}$ be a subsequential limit of $\{\bar{\mu}^n\}_{n\in\mathbb{N}^*}$. We have for all $t \geq 0$

$$\chi_t = \bar{\mu}_0^* + \int_0^t (\chi_s)' \, ds + \lambda t \bar{\alpha}^*,$$

which is nothing but the integrated transport equation associated to $(\bar{\mu}_0^*, g)$, where $g_t := \lambda t \bar{\alpha}^*$, $t \geq 0$. From Theorem 1, we have for all $t \geq 0$ and $\phi \in \mathcal{S}$

$$\langle \chi_t, \phi \rangle = \langle \tau_t \bar{\mu}_0^*, \phi \rangle + \langle g_t, \phi \rangle + \langle N(g_t), \phi \rangle$$
$$= \langle \bar{\mu}_0^*, \tau_t \phi \rangle + \lambda t \langle \bar{\alpha}^*, \phi \rangle - \lambda \int_0^t s \frac{d}{ds} \langle \bar{\alpha}^*, \tau_{t-s}\phi \rangle \, ds = \langle \bar{\mu}_t^*, \phi \rangle,$$

integrating by parts. The subsequential limit is therefore unique in the space $\mathcal{D}([0,\infty), \mathcal{M}_f^+)$, and equal to $\bar{\mu}^*$, since $\mathcal{S}$ is a separating class of $\mathcal{M}_f^+$. □

**5. Fluid limits of some performance processes.** Let us provide some applications of Theorem 2 to the asymptotic estimation of some performance processes describing the queueing system. Assume that Hypothesis 1 hold and, in addition, that the limiting initial profile $\bar{\mu}_0^*$ in Hypothesis 1 is such that

(27) $\qquad\qquad\qquad\qquad \bar{\mu}_0^*$ has no atom.

Assumption (27) implies in view of (20) that $\bar{\mu}_t^*$ has no atom for any $t \geq 0$. Thus, Theorem 2 and the Continuous Mapping Theorem yield that, for any $x < y$, the sequence $\{(\langle \bar{\mu}_t^n, \mathbf{1}_{(x,y)}\rangle)_{t\geq 0}\}_{n\in\mathbb{N}^*}$ tends in distribution to the real deterministic function $(\langle \bar{\mu}_t^*, \mathbf{1}_{(x,y)} \rangle)_{t \geq 0}$ (the boundaries may as well be included or infinite). Thus, a fluid limit approximation of the process counting the customers having residual service times in the range $(x, y)$ is given for all $t \geq 0$ by

$$\langle \bar{\mu}_t^*, \mathbf{1}_{(x,y)} \rangle = \bar{\mu}_0^*((x+t, y+t)) + \lambda \int_0^t \left( \int_{x+s}^{y+s} d\bar{\alpha}^*(x) \right) ds.$$

In particular, a fluid approximation of the normalized congestion process $\bar{X}^n$ is given by the process $\bar{X}^*$ defined for all $t$ by

$$\bar{X}_t^* = \langle \bar{\mu}_t^*, \mathbf{1}_{\mathbb{R}_+^*} \rangle = \bar{\mu}_0^*((t, \infty)) + \lambda \int_0^t \left( \int_s^{+\infty} d\bar{\alpha}^*(x) \right) ds,$$

whereas the normalized service process $\bar{S}^n$ can be approximated by $\bar{S}^*$, where

$$\bar{S}_t^* = \langle \bar{\mu}_t^*, \mathbf{1}_{\mathbb{R}^-} \rangle = \bar{\mu}_0^*((-\infty, t]) + \lambda \int_0^t \left( \int_0^s d\bar{\alpha}^*(x) \right) ds.$$

Remark now that for all $t \leq T$, all $x \geq 0$ and all $n \in \mathbb{N}^*$,

$$\bar{W}_t^n \leq \langle \bar{\mu}_t^n, I\mathbf{1}_{(0,x]} \rangle + \langle \bar{\mu}_0^n, I\mathbf{1}_{(t+x,\infty)} \rangle + \frac{1}{n} \sum_{i=1}^{N_{nt}^n} \bar{\sigma}_i^n \mathbf{1}_{\bar{\sigma}_i^n \geq x}.$$

In view of (24), the third process on the right tends to $(\lambda t \langle \bar{\alpha}^*, I\mathbf{1}_{(x,\infty)} \rangle)_{t \geq 0}$. Hence, it is bounded over compact sets and uniformly in $n$ since $\langle \bar{\alpha}^*, I\mathbf{1}_{(x,\infty)} \rangle$ is finite. The same argument applies to the second process on the right since $\langle \bar{\mu}_0^n, I\mathbf{1}_{(x,\infty)} \rangle$ is finite as well, and the first one tends in distribution to $(\langle \bar{\mu}_t^*, I\mathbf{1}_{(0,x]} \rangle)_{t \geq 0}$. Hence, $\{\bar{W}^n\}_{n \in \mathbb{N}^*}$ is tight, and any subsequential limit thus reads $(\langle \bar{\mu}_t^*, I\mathbf{1}_{\mathbb{R}_+^*} \rangle)_{t \geq 0}$. In other words, the fluid limit of $\bar{W}^n$ is given by $\bar{W}^*$, where

$$\bar{W}_t^* = \langle \bar{\mu}_t^*, I\mathbf{1}_{\mathbb{R}_+^*} \rangle$$
$$= \int_t^{+\infty} (x-t) \, d\bar{\mu}_0^*(x) + \lambda \int_0^t \left( \int_s^{+\infty} (x-s) \, d\bar{\alpha}^*(x) \right) ds.$$

## 6. Functional central limit theorem.

6.1. *Preliminary results.* In this section we prove two technical results (Propositions 3 and 4), which will be useful in the sequel. We refer again the reader to the definitions and notation introduced in Section 2. Throughout this whole section, fix $T > 0$.

LEMMA 1.  *For all $\phi \in \mathcal{S}$, the collection $\{\tau_r \phi, r \in [0, T]\}$ is a pre-compact set of $\mathcal{S}$.*

PROOF. For all semi-norm associated to the integers $\beta$, $\gamma$, for all $r \in [0, T]$,

$$|\tau_r \phi|_{\beta,\gamma} = \sup_x |x^\beta \phi^{(\gamma)}(x-r)| \leq \sum_{j=1}^\beta C_\beta^j r^{\beta-j} |\phi|_{j,\gamma} \leq \sum_{j=1}^\beta C_\beta^j T^{\beta-j} |\phi|_{j,\gamma}.$$

The set $\{\tau_r \phi, r \in [0,T]\}$ is thus bounded in the nuclear space $\mathcal{S}$: it is a pre-compact set in view of Proposition 4.4.7, page 81 of [23]. □

LEMMA 2. *For all $\phi \in \mathcal{S}$, the mapping $r \longmapsto \tau_r \phi$ is continuous from $\mathbb{R}$ into $\mathcal{S}$.*

PROOF. Let $r' < r$. Let $\beta$ and $\gamma \in \mathbb{N}$. For some $s \in ]r', r[$,

$$|\tau_{r'} \phi - \tau_r \phi|_{\beta,\gamma} = \sup_x |x^\beta (\phi^{(\gamma)}(x-r') - \phi^{(\gamma)}(x-r))|$$

$$\leq |r - r'| \sum_{j=1}^\beta C_\beta^j s^\beta |\phi|_{j,\gamma+1}$$

$$\leq |r - r'| \sum_{j=1}^\beta C_\beta^j (r+1)^\beta |\phi|_{j,\gamma+1} =: M_{r,\beta,\gamma,\phi}.$$

Thus, for all $\varepsilon > 0$, fixing $\eta := \varepsilon (M_{r,\beta,\gamma,\phi})^{-1}$, we have $|\tau_{r'} \phi - \tau_r \phi|_{\beta,\gamma} < \varepsilon$ whenever $|r - r'| < \eta$. □

We therefore have the following results.

PROPOSITION 3. *The mapping $G$ defined by (1) is continuous from $\mathcal{C}([0,T], \mathcal{S}')$ into itself.*

PROOF. We first prove that $G$ takes its values in $\mathcal{C}([0,T], \mathcal{S}')$. Let $\phi \in \mathcal{S}$, $X \in \mathcal{C}([0,T], \mathcal{S}')$, and let $G_\phi^X$ be the mapping from $\mathbb{R}^2$ into $\mathbb{R}$ defined by $G_\phi^X(t, r) = \langle X_t, \tau_r \phi \rangle$. Let $t, r \leq T$. For all $t', r'$,

$$(28) \quad |G_\phi^X(t', r') - G_\phi^X(t, r)| \leq |\langle X_{t'} - X_t, \tau_r \phi \rangle| + |\langle X_t, \tau_{r'} \phi - \tau_r \phi \rangle|.$$

On the one hand, $(\langle X_s, \phi \rangle)_{s \geq 0} \in \mathcal{D}([0,T], \mathbb{R})$, and is thus bounded on $[0, T]$:

$$\sup_{s \leq T} |\langle X_s, \phi \rangle| = \sup_{s \in [0,T] \cap \mathbb{Q}} |\langle X_s, \phi \rangle| < \infty.$$

Therefore, $\mathcal{X} := \{X_s, s \in [0,T] \cap \mathbb{Q}\}$ is a weakly bounded subset of the set of continuous mappings from $\mathcal{S}$ into $\mathbb{R}$. The space $\mathcal{S}$ is Fréchet, and hence is a tonnel in view of corollary 0, page III.25 of [4]. From the Banach–Steinhaus theorem (Theorem 1, page III.25 of [4]), $\mathcal{X}$ is equicontinuous. Hence, for all $\varepsilon > 0$,

$$V_\varepsilon := \bigcap_{s \in [0,T] \cap \mathbb{Q}} X_s^{-1}\left(\left]-\frac{\varepsilon}{2}, \frac{\varepsilon}{2}\right[\right) = \left\{\phi, \sup_{s \in [0,T] \cap \mathbb{Q}} |\langle X_s, \phi \rangle| < \frac{\varepsilon}{2}\right\}$$

is a neighborhood of $\mathbf{0}$, the zero of $\mathcal{S}$. The mapping $r \mapsto \tau_r \phi$ being continuous in view of Lemma 2, for some $\eta_{r,\phi} > 0$

$$|r' - r| < \eta_{r,\phi} \Longrightarrow (\tau_{r'}\phi - \tau_r\phi) \in V_\varepsilon \Longrightarrow |\langle X_t, \tau_{r'}\phi - \tau_r\phi \rangle| < \frac{\varepsilon}{2}.$$

On the other hand, $(\langle X_s, \tau_r\phi \rangle)_{s \geq 0} \in \mathcal{C}([0,T], \mathbb{R})$, hence, for some $\delta_{r,t,\phi}$ and all $t'$,

$$|t' - t| < \delta_{r,t,\phi} \Longrightarrow |\langle X_{t'} - X_t, \tau_r\phi \rangle| < \frac{\varepsilon}{2}.$$

For all $t, r$ and all $\varepsilon > 0$, the two previous relations in (28) imply that, for some $\eta_{r,\phi}$ and $\delta_{r,t,\phi}$,

$$|t' - t| < \delta_{r,t,\phi}, |r' - r| < \eta_{r,\phi} \Longrightarrow |G_\phi^X(t', r') - G_\phi^X(t, r)| < \varepsilon.$$

The mapping $G_\phi^X$ is thus continuous, and hence uniformly continuous on the compact set $([0,T])^2$: for some $\delta_\phi, \eta_\phi > 0$,

$$\sup_{r,r',|r-r'|<\eta_\phi} \sup_{t,t'|t-t'|<\delta_\phi} |G_\phi^X(t', r') - G_\phi^X(t, r)| < \varepsilon.$$

Finally, letting $\xi_\phi := \eta_\phi \wedge \delta_\phi$,

$$\sup_{t,t',|t-t'|<\xi_\phi} |\langle G(X)_{t'}, \phi \rangle - \langle G(X)_t, \phi \rangle|$$
$$\leq \sup_{r,r',|r-r'|<\eta_\phi} \sup_{t,t',|t-t'|<\delta_\phi} |G_\phi^X(t', r') - G_\phi^X(t, r)| < \varepsilon.$$

Thus, $G(X) \in \mathcal{C}([0,T], \mathcal{S}')$ since the map $t \mapsto \langle G(X)_t, \phi \rangle$ is continuous for all $\phi \in \mathcal{S}$.

Let us now show the continuity of $G$. Let $\{X^n\}_{n \in \mathbb{N}^*}$ be a sequence of $\mathcal{D}([0,T], \mathcal{S}')$ tending to $X$. In particular, for all $\phi \in \mathcal{S}$, $\{(\langle X_t^n, \phi \rangle)_{t \geq 0}\}_{n \in \mathbb{N}^*}$ tends in $\mathcal{D}([0,T], \mathbb{R})$ to $(\langle X_t, \phi \rangle)_{t \geq 0}$, and thus $\sup_{n \in \mathbb{N}^*} \sup_{t \leq T} |\langle X_t^n, \phi \rangle| < \infty$. Hence, in view of the Banach–Steinhaus theorem,

$$\mathcal{H} := \{X_t^n, n \in \mathbb{N}, t \in \mathbb{Q} \cap [0,T]\}$$

is equicontinuous. Therefore, for all $\phi \in \mathcal{S}$,

$$\sup_{t \leq T} |\langle X_t^n - X_t, \phi \rangle| \xrightarrow[n \to \infty]{} 0,$$

and consequently, for all $k$, all collection $\{\phi_i, i = 1, \ldots, k\}$ of $\mathcal{S}$ and all $i \leq k$,

$$\sup_{t \leq T} \max_{i=1,\ldots,k} |\langle X_t^n - X_t, \phi_i \rangle| \xrightarrow[n \to \infty]{} 0.$$

The latter means that for all semi-norm $p_w$ of the weak topology,

$$\sup_{t \leq T} p_w(X_t^n - X_t) \xrightarrow[n \to \infty]{} 0.$$

Thus, in view of Pietsch's theorem ([23], Proposition 0.6.7, page 9), for all pre-compact set $K$ of $\mathcal{S}$,

$$\sup_{t \leq T} \sup_{\phi \in K} |\langle X_t^n - X_t, \phi \rangle| \xrightarrow[n \to \infty]{} 0.$$

This result can be applied for a fixed $\phi \in \mathcal{S}$ to the collection $\{\tau_r \phi, r \in [0, T]\}$ (which is a pre-compact set of $\mathcal{S}$ from Lemma 1) to obtain

$$\sup_{t \leq T} \sup_{r \leq T} |\langle X_t^n - X_t, \tau_r \phi \rangle| \xrightarrow[n \to \infty]{} 0.$$

It follows that $G$ is continuous, since

$$\sup_{t \leq T} |\langle G(X^n)_t - G(X)_t, \phi \rangle| = \sup_{t \leq T} |\langle X_t^n - X_t, \tau_t \phi \rangle|$$

$$\leq \sup_{t \leq T} \sup_{r \leq T} |\langle X_t^n - X_t, \tau_r \phi \rangle| \xrightarrow[n \to \infty]{} 0. \quad \square$$

PROPOSITION 4. *The mapping $N$ defined by (3) is continuous from $\mathcal{D}([0, T], \mathcal{S}')$ into $\mathcal{C}([0, T], \mathcal{S}')$.*

PROOF. That $N$ takes values in $\mathcal{C}([0, T], \mathcal{S}')$ can be shown similarly to Proposition 3. For the continuity of $N$, remark that the mapping $H$ defined by (2) is continuous from $\mathcal{D}([0, T], \mathcal{S}')$ into $\mathcal{C}([0, T], \mathcal{S}')$ since, as well known, the mapping $(X_t)_{t \geq 0} \longmapsto (\int_0^t X_s \, ds)_{t \geq 0}$ is continuous from $\mathcal{D}([0, T], \mathbb{R})$ into $\mathcal{C}([0, T], \mathbb{R})$ for the Skorokhod topology. The proof then proceeds as that of Proposition 3. $\square$

6.2. *Central limit theorem.* Fix $\xi, \{\lambda^n\}_{n \in \mathbb{N}^*}, \lambda$, and $\{\sigma^n\}_{n \in \mathbb{N}^*}$ satisfying Hypothesis 1. $\{\bar{\mu}^n\}_{n \in \mathbb{N}^*}$ is the corresponding sequence of normalized profile processes. From Theorem 2, $\{\bar{\mu}^n\}_{n \in \mathbb{N}^*}$ tends in distribution to $\bar{\mu}^*$ in $\mathcal{D}([0, \infty), \mathcal{M}_f^+)$, defined by (20). Throughout this section, we make the following additional assumptions.

HYPOTHESIS 2.

(29) $$\sqrt{n}(\lambda^n - \lambda) \xrightarrow[n \to \infty]{} 0.$$

There exists $\mathcal{Y}_0 \in \mathcal{S}'$, such that for all $\phi \in \mathcal{S}$,

(30) $$\sqrt{n}\langle \bar{\mu}_0^n - \bar{\mu}_0^*, \phi \rangle \longrightarrow \langle \mathcal{Y}_0, \phi \rangle \quad \text{in probability},$$

(31) $$\sqrt{n}|\langle \bar{\alpha}^n, I \rangle - \langle \bar{\alpha}^*, I \rangle| \underset{n \to \infty}{\longrightarrow} 0.$$

Let $n \in \mathbb{N}^*$. Define the following processes for all $t$ by

$$\mathcal{Y}_t^n = \sqrt{n}(\bar{\mu}_t^n - \bar{\mu}_t^*),$$
$$\mathcal{M}_t^n = \sqrt{n}\bar{M}_t^n,$$
$$\mathcal{A}_t^n(\phi, \psi) = \lambda^n t \langle \bar{\alpha}^n, \phi\psi \rangle, \qquad \phi, \psi \in \mathcal{S}.$$

The process $\mathcal{M}^n$ is an $\mathcal{S}'$-valued $\mathcal{G}_t^n$-martingale, hence, for all $\phi \in \mathcal{S}$, the process $(\langle \mathcal{M}_t^n, \phi \rangle)_{t \geq 0}$ is a real square integrable martingale of increasing process $(\mathcal{A}_t^n(\phi, \phi))_{t \geq 0}$ [by (19)].

LEMMA 3. *Under Hypotheses 1 and 2, for all $T > 0$, the sequence $\{\mathcal{M}^n\}_{n \in \mathbb{N}^*}$ converges in distribution to the martingale $\mathcal{M}$ of $\mathcal{D}([0,T], \mathcal{S}')$, whose tensor-quadratic process satisfies for all $t \in [0,T]$ and all $\phi, \psi \in \mathcal{S}$*

(32) $$<< \mathcal{M} >>_t (\phi, \psi) = \lambda t \langle \bar{\alpha}^*, \phi\psi \rangle.$$

PROOF. We use the convergence criterion, Theorem A.3 in the Appendix. Define for all $t \geq 0$ and $\phi, \psi \in \mathcal{S}$

$$\gamma_t(\phi, \psi) := \lambda t \langle \bar{\alpha}^*, \phi\psi \rangle.$$

For all $\phi \in \mathcal{S}$, $\gamma_0(\phi, \phi) = 0$ and for all $n$, $((\langle \mathcal{M}_t^n, \phi \rangle)^2 - \mathcal{A}_t^n(\phi, \phi))_{t \geq 0}$ is a $\mathcal{G}_t^n$-local martingale, thus conditions (37) and (38) are verified. On the other hand, for all $T > 0$ and $\psi \in \mathcal{S}$, we have **P**-a.s. for all $t \leq T$

$$\mathcal{A}_t^n(\phi, \psi) - \gamma_t(\phi, \psi) = (\lambda^n - \lambda) t \langle \bar{\alpha}^n, \phi\psi \rangle + \lambda t \langle \bar{\alpha}^n - \bar{\alpha}^*, \phi\psi \rangle,$$

hence,

$$\sup_{t \leq T} \{\mathcal{A}_t^n(\phi, \psi) - \gamma_t(\phi, \psi)\}^2$$
$$\leq 2T^2 \{(\lambda^n - \lambda)^2 \|\phi\psi\|_\infty^2 + \lambda^2 \|(\phi\psi)'\|_\infty^2 \langle \bar{\alpha}^n - \bar{\alpha}^*, I \rangle\} \underset{n \to \infty}{\longrightarrow} 0,$$

which proves (41) taking $\psi := \phi$. Finally, for all $T > 0$, condition (40) is met since $(\mathcal{A}_t^n(\phi, \phi))_{t \geq 0} \in \mathcal{C}([0,T], \mathbb{R})$ and

$$\mathbf{E}\left[\sup_{t \leq T}(\langle \mathcal{M}_t^n, \phi \rangle - \langle \mathcal{M}_{t-}^n, \phi \rangle)^2\right] \mathbf{E}\left[\sup_{t \leq T}(\langle \bar{\mu}_t^n, \phi \rangle - \langle \bar{\mu}_{t-}^n, \phi \rangle)^2\right]$$
$$\leq \frac{\|\phi\|_\infty}{n} \underset{n \to \infty}{\longrightarrow} 0,$$

the jumps of $(\langle \bar{\mu}_t^n, \phi \rangle)_{t \geq 0}$ being a.s. of size less than $\frac{\|\phi\|_\infty}{n}$: we have (39) as well. Therefore, we can apply Theorem A.3 to state that

(33) $$\{(\langle \mathcal{M}_t^n, \phi \rangle)_{t \geq 0}\}_{n \in \mathbb{N}^*} \Longrightarrow (P_t)_{t \geq 0} \quad \text{in } \mathcal{D}([0, \infty), \mathbb{R}),$$

where $P$ is a continuous martingale of increasing process $(\gamma_t(\phi, \phi))_{t \geq 0}$. In particular, for all $T > 0$ and all $\phi \in \mathcal{S}$, $\{(\langle \mathcal{M}_t^n, \phi \rangle)_{t \geq 0}\}_{n \in \mathbb{N}^*}$ is tight in $\mathcal{D}([0, T], \mathbb{R})$, and hence (Theorem A.1 in the Appendix), $\{\mathcal{M}^n\}_{n \in \mathbb{N}^*}$ is tight in $\mathcal{D}([0, T], \mathcal{S}')$. All subsequential limit $\mathcal{M}$ thus satisfies $\langle \mathcal{M}_t, \phi \rangle = P_t$ for all $t \in [0, T]$ and all $\phi \in \mathcal{S}$, which means with (33) that $\{\mathcal{M}^n\}_{n \in \mathbb{N}^*}$ tends in distribution to the martingale $\mathcal{M}$ of $\mathcal{D}([0, T], \mathcal{S}')$. Finally and according to the previous arguments,

$$\{(<\langle \mathcal{M}_\cdot^n, \phi \rangle, \langle \mathcal{M}_\cdot^n, \psi \rangle >_t)_{t \geq 0}\}_{n \in \mathbb{N}^*} = \{(\mathcal{A}_t^n(\phi, \psi))_{t \geq 0}\}_{n \in \mathbb{N}^*}$$

tends to $(<< \mathcal{M} >>_t (\phi, \psi))_{t \geq 0}$ as well as to $(\gamma_t(\phi, \psi))_{t \geq 0}$, which completes the proof. □

Consider now the Hilbert space $\mathbf{L}^2(d\bar{\alpha}^*)$. Since $\bar{\alpha}^*$ is a probability measure, we have $\mathcal{S} \subset \mathbf{L}^2(d\bar{\alpha}^*)$, and hence, $(\mathbf{L}^2(d\bar{\alpha}^*))' \subset \mathcal{S}'$. Moreover, $\mathbf{L}^2(d\bar{\alpha}^*)$ admits a countable basis, denoted $\{h_i\}_{i \in \mathbb{N}}$. For all $i \geq 0$ and all $\psi \in \mathbf{L}^2(d\bar{\alpha}^*)$, denote

$$c_i(\psi) = \int \psi(x) h_i(x) \, d\bar{\alpha}^*(x),$$

the $i$th coordinate of $\psi$ in that basis. We have the following result.

THEOREM 3. *Under Hypotheses 1 and 2, for all $T > 0$, the sequence $\{\mathcal{Y}^n\}_{n \in \mathbb{N}^*}$ tends in distribution in $\mathcal{D}([0, T], \mathcal{S}')$ to the process $\mathcal{Y}^*$ defined for all $t \geq 0$ and all $\phi \in \mathcal{S}$ by*

(34) $$\langle \mathcal{Y}_t^*, \phi \rangle = \langle \mathcal{Y}_0, \tau_t \phi \rangle + \sqrt{\lambda} \sum_{i \geq 0} \int_0^t c_i(\tau_{t-s} \phi) \, dB_s^i,$$

*where $\{B^i\}_{i \geq 0}$ is a sequence of independent real standard Brownian motions.*

PROOF. For all $n \in \mathbb{N}$ and all $t \in [0, T]$, **P**-a.s.,

$$\mathcal{Y}_t^n = \mathcal{Y}_0^n + \int_0^t (\mathcal{Y}_s^n)' \, ds + \sqrt{n}\{\lambda^n t \bar{\alpha}^n - \lambda t \bar{\alpha}^*\} + \mathcal{M}_t^n.$$

Hence, $(\mathcal{Y}_t^n)_{t \geq 0}$ solves a transport equation, and in view of Theorem 1, the above equation amounts **P**-a.s. for all $t \leq T$ to

$$\mathcal{Y}_t^n = G(\mathcal{Y}_0^n)_t + \sqrt{n}\{\lambda^n t \bar{\alpha}^n - \lambda t \bar{\alpha}^*\} + \mathcal{M}_t^n$$
$$\quad + \sqrt{n} \lambda^n N(\cdot \bar{\alpha}^n) - \sqrt{n} \lambda N(\cdot \bar{\alpha}^*) + N(\mathcal{M}^n)_t$$
(35) $$= G(\mathcal{Y}_0^n)_t + \sqrt{n}(\lambda^n - \lambda) t \bar{\alpha}^n - \lambda t \sqrt{n}(\bar{\alpha}^* - \bar{\alpha}^n) + \mathcal{M}_t^n$$
$$\quad + \sqrt{n}(\lambda^n - \lambda) N(\cdot \bar{\alpha}^n)_t - \lambda\{N(\sqrt{n} \cdot \bar{\alpha}^*)_t - N(\sqrt{n} \cdot \bar{\alpha}^n)_t\}$$
$$\quad + N(\mathcal{M}^n)_t,$$

where $(\mathcal{Y}_0^n)_{t\geq 0}$ denotes the $\mathcal{S}'$-valued process constantly equal to $\mathcal{Y}_0^n$ and $\cdot\bar{\alpha}^n$ (resp. $\cdot\bar{\alpha}^*$, $\sqrt{n}\cdot\bar{\alpha}^n$, $\sqrt{n}\cdot\bar{\alpha}^*$) denotes the $\mathcal{S}'$-valued process $(t\bar{\alpha}^n)_{t\geq 0}$ [resp. $(t\bar{\alpha}^*)_{t\geq 0}$, $(\sqrt{n}t\bar{\alpha}^n)_{t\geq 0}$, $(\sqrt{n}t\bar{\alpha}^*)_{t\geq 0}$]. First, from (30), for all $\phi \in \mathcal{S}$, $\xi > 0$,

$$\mathbf{P}[|\langle \mathcal{Y}_0^n - \mathcal{Y}_0, \phi\rangle| \geq \xi] \xrightarrow[n\to\infty]{} 0.$$

Hence, in view of Corollary A.1 in the Appendix, we have that $(\mathcal{Y}_0^n)_{t\geq 0} \Rightarrow (\mathcal{Y}_0)_{t\geq 0}$ in $\mathcal{D}([0,T], \mathcal{S}')$. Thus, from Proposition 3 and the Continuous Mapping Theorem (see [2]),

$$(G(\mathcal{Y}_0^n)_t)_{t\geq 0} \Rightarrow (G(\mathcal{Y}_0)_t)_{t\geq 0} \quad \text{in } \mathcal{D}([0,T], \mathcal{S}').$$

On the other hand, we have

$$\sup_{0\leq t\leq T} \sqrt{n}|\lambda^n - \lambda|t|\langle\bar{\alpha}^n,\phi\rangle| \leq \sqrt{n}|\lambda^n - \lambda|T\|\phi\|_\infty \xrightarrow[n\to\infty]{} 0,$$

$$\sup_{0\leq t\leq T} \sqrt{n}|\lambda^n - \lambda||\langle N(\cdot\bar{\alpha}^n)_t,\phi'\rangle| \leq \sqrt{n}|\lambda^n - \lambda|T^2\|\phi'\|_\infty \xrightarrow[n\to\infty]{} 0,$$

$$\sup_{0\leq t\leq T} \lambda t\sqrt{n}|\langle\bar{\alpha}^*,\phi\rangle - \langle\bar{\alpha}^n,\phi\rangle| \leq \lambda T\sqrt{n}\|\phi'\|_\infty|\langle\bar{\alpha}^*,I\rangle - \langle\bar{\alpha}^n,I\rangle| \xrightarrow[n\to\infty]{} 0.$$

Hence, the sequences $\{(\lambda^n - \lambda)\sqrt{n}\cdot\bar{\alpha}^n\}_{n\in\mathbb{N}^*}$, $\{(\lambda^n - \lambda)N(\sqrt{n}\cdot\bar{\alpha}^n)\}_{n\in\mathbb{N}^*}$ and $\{\lambda(\sqrt{n}\cdot\bar{\alpha}^* - \sqrt{n}\cdot\bar{\alpha}^n)\}_{n\in\mathbb{N}^*}$ converge to the null process $(\mathbf{0})_{t\geq 0}$ of $\mathcal{C}([0,T], \mathcal{S}')$. The convergence of this last sequence implies in view of Proposition 4 that

$$\lambda\{N(\sqrt{n}\cdot\bar{\alpha}^*) - N(\sqrt{n}\cdot\bar{\alpha}^n)\} \xrightarrow[n\to\infty]{} (\mathbf{0})_{t\geq 0} \quad \text{in } \mathcal{C}([0,T], \mathcal{S}').$$

Now, from Lemma 3, $\mathcal{M}^n \Longrightarrow \mathcal{M}$ in $\mathcal{D}([0,T], \mathcal{S}')$, thus, in view of Proposition 4 and the Continuous Mapping Theorem, $N(\mathcal{M}^n) \Longrightarrow N(\mathcal{M})$ in $\mathcal{D}([0,T], \mathcal{S}')$. Consequently, $\{\mathcal{Y}^n\}_{n\in\mathbb{N}^*}$ is tight and from (35), its limit in distribution $\mathcal{Y}^*$ satisfies $\mathbf{P}$-a.s. for all $t \in [0,T]$:

$$(36) \quad \mathcal{Y}_t^* = G(\mathcal{Y}_0)_t + \mathcal{M}_t + N(\mathcal{M})_t = \tau_t\mathcal{Y}_0 + \mathcal{M}_t + \int_0^t \tau_{t-s}(\mathcal{M}_s)'\,ds.$$

In view of (32), the process $\mathcal{M}$ reads $\mathcal{M} = \sqrt{\lambda}\mathbb{B}(\bar{\alpha}^*)$, where $\mathbb{B}(\bar{\alpha}^*)$ is the cylindrical Brownian motion on $(\mathbf{L}^2(d\bar{\alpha}^*))'$ (see [18]). Thus, from (36), we have for all $\phi \in \mathcal{S}$ and for all $t$

$$\langle \mathcal{Y}_t^*, \phi\rangle = \langle \mathcal{Y}_0, \tau_t\phi\rangle + \sqrt{\lambda}\bigg\{\langle \mathbb{B}(\bar{\alpha}^*)_t, \phi\rangle - \int_0^t \langle \mathbb{B}(\bar{\alpha}^*)_s, \tau_{t-s}\phi'\rangle\,ds\bigg\}$$

$$= \langle \mathcal{Y}_0, \tau_t\phi\rangle + \sqrt{\lambda}\sum_{i\geq 0}\bigg\{c_i(\phi)B_t^i - \int_0^t c_i(\tau_{t-s}\phi')B_s^i\,ds\bigg\},$$

where $\{B^i\}_{i\geq 0}$ is a sequence of independent real standard Brownian motions (see again [18]). Note that the latter series converges in $\mathbf{L}^2(\Omega)$ since $\phi \in \mathbf{L}^2(d\bar{\alpha}^*)$. Remark now that, for all $i \geq 0$,

$$\frac{d(c_i(\tau_{t-s}\phi))}{ds} = c_i(\tau_{t-s}\phi').$$

Using Itô's integration by parts formula, this yields to

$$c_i(\phi)B_t^i - \int_0^t c_i(\tau_{t-s}\phi')B_s^i\,ds = \int_0^t c_i(\tau_{t-s}\phi)\,dB_s^i,$$

which completes the proof. □

6.3. *Diffusion approximations of the performance processes.* In this section we show that Theorem 3 can be used to provide diffusion approximations for the congestion, service and workload processes. The stochastic integrals on the right-hand side of (34) make sense for any $\phi \in \mathbf{L}^2(d\bar{\alpha}^*)$, hence, the process $(\mathcal{Y}_t^*)_{t\geq 0}$ takes values in $(\mathbf{L}^2(d\bar{\alpha}^*))'$, provided that the limiting initial state $\mathcal{Y}_0$ belongs to $(\mathbf{L}^2(d\bar{\alpha}^*))'$. Showing that in that case the convergence announced in Theorem 3 holds in $\mathcal{D}([0,T],(\mathbf{L}^2(d\bar{\alpha}^*))')$ is an open problem at this point, that is beyond the scope of this paper. However, it is easily seen that under Hypothesis 1 the limiting service time distribution $\bar{\alpha}^*$ is such that the functions $\mathbf{1}_{\mathbb{R}_+^*}$ and $\mathbf{1}_{\mathbb{R}_-}$ belong to $\mathbf{L}^2(d\bar{\alpha}^*)$. Hence, diffusion approximations for the congestion and service processes can be provided by the real processes obtained when fixing $\phi = \mathbf{1}_{\mathbb{R}_+^*}$ and $\phi = \mathbf{1}_{\mathbb{R}_-}$ in (34), that is, for all $t \geq 0$,

$$\langle \mathcal{Y}_t^*, \mathbf{1}_{\mathbb{R}_+^*} \rangle = \mathcal{Y}_0((t,+\infty)) + \sqrt{\lambda} \sum_{i\geq 0} \int_0^t \left( \int_{t-s}^{+\infty} h_i(x)\,d\bar{\alpha}^*(x) \right) dB_s^i$$

and

$$\langle \mathcal{Y}_t^*, \mathbf{1}_{\mathbb{R}_-} \rangle = \mathcal{Y}_0((-\infty,t]) + \sqrt{\lambda} \sum_{i\geq 0} \int_0^t \left( \int_0^{t-s} h_i(x)\,d\bar{\alpha}^*(x) \right) dB_s^i.$$

Moreover, provided that $\int_0^{+\infty} x^2\,d\bar{\alpha}^*(x) < \infty$ (i.e., the limiting service time has a finite second moment), a diffusion approximation for the workload process is given for all $t \geq 0$ by

$$\langle \mathcal{Y}_t^*, I\mathbf{1}_{\mathbb{R}_+^*} \rangle = \int_t^{+\infty} (x-t)\,d\mathcal{Y}_0(x)$$
$$+ \sqrt{\lambda} \sum_{i\geq 0} \int_0^t \left( \int_{t-s}^{+\infty} (x+s-t)h_i(x)\,d\bar{\alpha}^*(x) \right) dB_s^i.$$

Assume, for example, that $\bar{\alpha}^*$ is the exponential distribution $\varepsilon(1)$ [which is the limiting distribution obtained when taking $\sigma^n \sim \varepsilon(1/n)$ for all $n \in \mathbb{N}^*$]. In that case, it is well known that a possible $\{h_i\}_{i\in\mathbb{N}}$ is given by the sequence of Laguerre's polynomials, defined for all $i \geq 0$ by

$$h_i(x) = \frac{e^x}{i!} \frac{d^i}{dx^i}(e^{-x}x^i).$$

Then a diffusion approximation of the congestion process (and accordingly, of the service and workload processes) is given for all $t \geq 0$ by

$$\langle \mathcal{Y}_t^*, \mathbf{1}_{\mathbb{R}_+^*} \rangle = \mathcal{Y}_0((t, +\infty)) + \sqrt{\lambda} \sum_{i \geq 0} \int_0^t \left( \int_{t-s}^{+\infty} \frac{1}{i!} \frac{d^i}{dx^i}(e^{-x} x^i) \, dx \right) dB_s^i.$$

## APPENDIX: TIGHTNESS AND WEAK CONVERGENCE ON FUNCTIONAL SPACES

Let us recall some results about tightness and weak convergence on metric spaces.

THEOREM A.1 (cf. [19], Theorem 4.1). *Let $F$ be a nuclear Fréchet space, $F'$ be its topological dual, and $\{X^n\}_{n \in \mathbb{N}^*}$ be a sequence of $\mathcal{D}([0,T], F')$-valued r.v. Then $\{X^n\}_{n \in \mathbb{N}^*}$ is tight if for all $\phi \in F$, the sequence $\{(\langle X_t^n, \phi \rangle)_{t \geq 0}\}_{n \in \mathbb{N}^*}$ is tight in $\mathcal{D}([0,T], \mathbb{R})$.*

This shows, in particular, since $\mathcal{S}$ is a nuclear Fréchet space, that the tightness of $\{X^n\}_{n \in \mathbb{N}^*}$ in $\mathcal{D}([0,T], \mathcal{S}')$ is granted by that of $\{(\langle X^n, \phi \rangle)_{t \geq 0}\}_{n \in \mathbb{N}^*}$ for all $\phi \in \mathcal{S}$. Moreover, since $\mathcal{S}$ is naturally a separating class of its topological dual $\mathcal{S}'$, we have the following:

COROLLARY A.1. *For all sequence $\{X^n\}_{n \in \mathbb{N}^*}$ of $\mathcal{D}([0,T], \mathcal{S}')$, $X^n \Longrightarrow X$ if, for all $\phi \in \mathcal{S}$, $(\langle X_t^n, \phi \rangle)_{t \geq 0} \Longrightarrow (\langle X_t, \phi \rangle)_{t \geq 0}$ in $\mathcal{D}([0,T], \mathcal{S}')$.*

The space $\mathcal{M}_f^+$ has no such good topological properties as that of $\mathcal{S}'$. Nevertheless, a tightness criterion is established for sequences of $\mathcal{M}_f^+$-valued processes.

THEOREM A.2 (Jakubowski's criterion, cf. [5], Theorem 3.6.4). *The sequence $\{X^n\}_{n \in \mathbb{N}^*}$ of $\mathcal{D}([0,T], \mathcal{M}_f^+)$ is tight if the following two conditions hold:*

**C.1.** *For all $\phi \in \mathcal{C}_b^1$, the sequence $\{(\langle X_t^n, \phi \rangle)_{t \geq 0}\}_{n \in \mathbb{N}^*}$ is tight in $\mathcal{D}([0, \infty), \mathbb{R})$,*
**C.2.** *For all $T > 0$ and $0 < \eta < 1$, there exists a compact subset $\mathbf{K}_{T,\eta}$ of $\mathcal{M}_f^+$ such that*

$$\liminf_{n \to \infty} \mathbf{P}[X_t^n \in \mathbf{K}_{T,\eta} \, \forall t \in [0,T]] \geq 1 - \eta.$$

Let us finally give the following result, which establishes a criterion for the convergence in distribution of a sequence of real processes to a diffusion process.

THEOREM A.3 ([10], Theorem 1.4). *Let $\{R^n\}$ and $\{A^n\}$ be two sequences of $\mathcal{D}([0,\infty),\mathbb{R})$, and an increasing function $\gamma \in \mathcal{D}([0,\infty),\mathbb{R})$ such that*

$$(37) \qquad R^n(0) = 0 \quad and \quad \gamma(0) = 0,$$

$$(38) \qquad R^n \text{ and } (R^n)^2 - A^n \text{ are two } \mathcal{F}_t - \text{local martingales},$$

*and for all $T > 0$,*

$$(39) \qquad \mathbf{E}\left[\sup_{t<T}\{R^n(t) - R^n(t-)\}^2\right] \underset{n\to\infty}{\longrightarrow} 0,$$

$$(40) \qquad \mathbf{E}\left[\sup_{t<T}\{A^n(t) - A^n(t-)\}\right] \underset{n\to\infty}{\longrightarrow} 0,$$

$$(41) \qquad \forall t > 0, \forall \varepsilon > 0 \qquad \mathbf{P}[\{A^n(t) - \gamma(t)\} > \varepsilon] \underset{n\to\infty}{\longrightarrow} 0.$$

*Then, $Y^n$ tends in distribution to a continuous martingale of increasing process $(\gamma(t))_{t\geq 0}$.*

GET/TÉLÉCOM PARIS—CNRS UMR 5141
46 RUE BARRAULT
75634 PARIS
FRANCE
E-MAIL: Laurent.Decreusefond@enst.fr

CEREMADE—CNRS UMR 5745
UNIVERSITÉ PARIS-DAUPHINE
PLACE DU MARÉCHAL DE LATTRE DE TASSIGNY
75016 PARIS
FRANCE
AND
LABORATOIRE DE MATHÉMATIQUES APPLIQUÉES
UNIVERSITÉ DE TECHNOLOGIE DE COMPIÈGNE
CENTRE DE RECHERCHES DE ROYALLIEU
BP 20 529
60 205 COMPIÈGNE CEDEX
FRANCE
E-MAIL: Pascal.Moyal@utc.fr